\documentclass[reqno]{amsart}
\usepackage[T1]{fontenc}
\usepackage[latin1]{inputenc}
\usepackage[english]{babel}
\usepackage[babel]{csquotes}

\usepackage{amssymb}
\usepackage{amsmath}
\usepackage{amsthm}
\usepackage{latexsym}
\usepackage{graphicx}
\usepackage{mathrsfs}
\usepackage{mathrsfs}
\usepackage{bbm}
\usepackage{verbatim}
\usepackage[usenames,dvipsnames]{color} 

\newtheorem{thm}{Theorem}[section]

\newtheorem{lem}[thm]{Lemma}
\newtheorem{proposition}[thm]{Proposition}

\numberwithin{equation}{section}

\def\enne{\mathbb{N}}
\def\erre{\mathbb{R}}
\def\Rz{\mathbb{R}}

\def\P{\mathbb{P}}

\def\E{\mathop{{}\mathbb{E}}}
\def\cL{\mathscr{L}}
\def\cF{\mathscr{F}}
\def\cB{\mathscr{B}}
\def\eps{\varepsilon}

\def\OO{\mathcal{O}}
\def\beq{\begin{equation}}
\def\eeq{\end{equation}}
\def\to{\rightarrow}

\def\embed{\hookrightarrow}

\def\norm #1{\left\|#1\right\|}

\def\sp #1#2{\left<#1,#2\right>}
\newcommand\ip\sp

\renewcommand{\d}{{\mathrm d}}

\begin{document}
\title[An order approach to SPDEs ]
{An order approach to SPDEs with antimonotone terms}

\author{Luca Scarpa}
\address[Luca Scarpa]{Faculty of Mathematics, University of Vienna, 
Oskar-Morgenstern-Platz 1, 1090 Wien, Austria.}
\email{luca.scarpa@univie.ac.at}
\urladdr{http://www.mat.univie.ac.at/$\sim$scarpa}

\author{Ulisse Stefanelli}
\address[Ulisse Stefanelli]{Faculty of Mathematics, University of Vienna, 
Oskar-Morgenstern-Platz 1, 1090 Wien, Austria  and  Istituto di Matematica
Applicata e Tecnologie Informatiche \textit{{E. Magenes}}, v. Ferrata 1, 27100
Pavia, Italy.}
\email{ulisse.stefanelli@univie.ac.at}
\urladdr{http://www.mat.univie.ac.at/$\sim$stefanelli}



\begin{abstract}
We consider a class of parabolic stochastic partial differential
equations featuring an antimonotone nonlinearity. The existence of
unique maximal and minimal
 variational  solutions is proved via a fixed-point argument for nondecreasing
mappings in ordered spaces. This relies on the validity of a
comparison principle.  
\end{abstract}

\maketitle

%%%%%%%%%%%%%%%%%%%%%%%%%%%%%%%%%%%%%%%%%%%%%%%%

\section{Introduction}
\setcounter{equation}{0}
\label{sec:intro}

This note is concerned with the existence of solutions to a class of 
parabolic stochastic partial differential equations (SPDEs).

The typical setting that we have in mind is the equation
\begin{equation}
  \label{eq:i0}
  \d u - {\rm div} (a(\nabla u)) \, \d t - b (u) \, \d t = f(u) \, \d t + G(u) \, \d W
  \qquad\text{in } (0,T)\times\OO
\end{equation}
suitably coupled with boundary and initial conditions,
with $\OO$ being a smooth bounded domain of $\erre^d$
and $T>0$ a fixed final time.

Here, the real-valued variable $u$ is 
defined on $\Omega\times[0,T]\times\OO$,
$a$ is monotone and polynomial, 
$f$ is Lipschitz continuous, and $G$ is a Lipschitz-type operator,
stochastically integrable with respect to $W$, a
cylindrical Wiener process on the underlying probability space
$(\Omega,\cF, \P)$. The function $b:\Rz \to \Rz$ is nondecreasing,
possibly being nonsmooth, so that the corresponding term in the
left-hand side of the equation is indeed antimonotone.

Our aim is to prove that a variational formulation of relation \eqref{eq:i0} admits a solution,
whenever complemented with suitable initial and boundary
conditions. If $b$ is Lipschitz continuous or $-b$ is nondecreasing
and continuous such existence follows from the
classical theory by {\sc Pardoux} \cite{Pard} and 
{\sc Krylov-Rozovski{\u\i}} \cite{KR-spde}, 
see also {\sc Liu--R\"ockner} \cite{LiuRo}. 
By contrast, we focus here in the case of $b$ linearly bounded but not 
continuous nor nondecreasing.

This situation, to the best of our
knowledge, has yet to be addressed. Indeed, 
the possible discontinuity of $-b$ prevents it
from being even locally Lipschitz-continuous,
hence also the refined well-posedness results
for SPDEs with locally monotone or locally 
Lipschitz-continuous drift (see again \cite{LiuRo}) cannot be applied.

The case of a nondecreasing but not Lipschitz continuous nonlinearity $b$
in \eqref{eq:i0} prevents from proving existence by a standard regularization or
approximation approach. In fact, the usual parabolic compactness seem to
be of little use in order to pass to the limit in the antimonotone
term $-b (u)$. We resort here in tackling the problem in an
ordered-space framework instead, by exploiting the fact that $b$ is
nondecreasing. 

At first, we check the validity of a comparison principle
by extending to the nonlinear frame of relation \eqref{eq:i0} the corresponding result
by {\sc Chekroun, Park,\& Temam}
\cite{chek-park-tem}, see Proposition~\ref{thm:1}. This comparison
principle allows us to reformulate the existence issue as a fixed-point problem for
nondecreasing mappings in ordered spaces. By implementing this
fixed-point procedure, we check in Theorem
\ref{thm:2} that equation \eqref{eq:i0} admits variational solutions.

The variational solutions that we obtain via such order method 
are considered in a {\em strong} probabilistic sense, i.e.~not
changing the original stochastic basis and Wiener process.
Let us stress that this 
is extremely satisfactory especially because 
no uniqueness is to be expected for the equation \eqref{eq:i0}.
Consequently, if one tackled the problem through 
classical approximation procedures and passage to the limit
by stochastic compactness arguments, 
the nonuniqueness of the limit problem would  
prevent from obtaining probabilistically strong solutions
by the classical procedure \`a la {\sc Yamada--Watanabe} \cite{yama-wata}.
The order argument that we employ is thus efficient in passing by
this problem and providing solutions in a strong probabilistic sense
even if no uniqueness is expected.
Still, one can prove that the set of solutions admits
unique maximal and minimal elements in the sense of the pointwise almost-everywhere order. 

Before going on, let us mention that order methods for
proving existence for SPDEs have already been used in the frame of
viscosity solvability. The reader is referred to the seminal papers by
{\sc Lions \& Souganidis} \cite{LS98,LS98b} as well to 
\cite{BM1,BM2} for a
collection of results
in this direction.
The novelty here is that we focus on weak solutions
instead and that comparison is combined with a fixed-point
procedure. The fixed-point Lemma \ref{kolo} corresponds indeed to an
abstract version of 
Perron's method.

The setting of the problem is discussed in Section \ref{sec:main} where
we collect some preliminaries and
we state our main results, namely Proposition \ref{thm:1}
(comparison principle) and Theorem  \ref{thm:1} (existence). The
corresponding proofs are then given in Sections
\ref{sec:thm1}-\ref{sec:thm2}, respectively.

%%%%%%%%%%%%%%%%%%%%%%%%%%%%%%%%%%%%%%%%%%%%%%%%

\section{Setting and main results}
\setcounter{equation}{0}
\label{sec:main}

The aim of this section is to specify assumptions and introduce a variational
formulation for equation 
\eqref{eq:i0}, by possibly allowing for additional dependencies in the
nonlinearities. Eventually, our main results 
Proposition \ref{thm:1} and Theorem \ref{thm:2}
are also presented.

Let $(\Omega, \cF, (\cF_t)_{t\in[0,T]}, \P)$ be a complete filtered probability space,
where $T>0$ is a given final time,  $W$ be a cylindrical Wiener process
on a separable Hilbert space $U$, and  
fix a complete orthonormal system $(e_k)_{k\in\enne}$ of $U$.
The {\it progressive} $\sigma$-algebra on $\Omega\times[0,T]$
is denoted by $\mathcal P$.
For any Banach space $E$ and $r,s\in\mathopen[1,\infty\mathclose)$
we denote by $L^r(\Omega; E)$ and $L^r(0,T; E)$ the usual 
functional spaces of Bochner $r$-integrable functions and
by $L^r_{\mathcal P}(\Omega; L^s(0,T; E))$ the space of
progressively measurable processes $\varphi:\Omega\times[0,T]\to E$
such that $$\E\left(\int_0^T\norm{\varphi(t)}_E^s\,\d t\right)^{r/s}<\infty\,.$$
For any pair of separable Hilbert spaces $E_1$ and $E_2$, the symbol
$\cL^2(E_1,E_2)$ denotes the space of Hilbert-Schmidt operators from 
$E_1$ to $E_2$.
 
Let $\OO\subset\erre^d$ be
nonempty, open, bounded set with Lipschitz boundary.
We define the separable Hilbert space
\begin{itemize}
\item[(S1)] $H:=L^2(\OO)$,
\end{itemize}
and endow it with its usual scalar product $(\cdot,\cdot)_H$ and norm
$\norm{\cdot}_H$. Moreover, we ask
\begin{itemize}
\item[(S2)] $V$ to be a separable reflexive Banach space, 
continuously and densely embedded in $H$, that $V$ and $V^*$ are
uniformly convex and that  $ V\embed L^4(\OO)$ continuously.
\end{itemize} 
Throughout the paper, we identify $H$ with its dual $H^*$ through its Riesz isomorphism,
so that the inclusions
\[
  V \embed H \embed V^* 
\]
are continuous and dense. The norm in $V$ and
the duality between $V^*$ and $V$
will be denoted by $\norm{\cdot}_V$ and $\ip{\cdot}{\cdot}_V$,
respectively. 

Assumption (S2) is fulfilled for each closed subspace
of $W^{1,p}(\OO)$ for $p \geq 4d /(4+d)$. In particular, homogeneous Dirichlet boundary
conditions could be complemented to \eqref{eq:i0} by letting 
$u\in V=W^{1,p}_0(\OO)$, other choices being obviously possible. The requirement on $V$ and $V^*$ being uniformly convex
relates to the validity of a suitable It\^o's formula,
\cite[Thm.~4.1--4.2]{Pard}. 

By allowing additional dependencies, we let the
nonlinear function $a: \Omega \times [0,T]\times \Rz^d \to \Rz^d$ in
\eqref{eq:i0} possibly depend on time and realization as well. In
particular, we ask $a$ to be a
Carath\'eodory function, monotone and with $p$-growth with respect to
the last variable. This allows to define the operator
$A:\Omega\times[0,T]\times V \to V^*$ as
$$\langle Au, v \rangle_V : = \int_\OO a(\omega,t,\nabla u ) \cdot \nabla v \, \d
x \quad \forall u,\, v \in V\,.$$
By referring now directly to such operator, we assume the following
\begin{itemize}
  \item[(A1)] $A:\Omega\times[0,T]\times V\to V^*$ is $\mathcal P\otimes\cB(V)$--$\cB(V^*)$ measurable;
  \item [(A2)] for every $(\omega,t)\in\Omega\times[0,T]$ and $\varphi,\psi,\zeta\in V$
  the map $r\mapsto\ip{A(\omega,t,\varphi + r\psi)}{\zeta}_V$, $r\in\erre$, is continuous;
  \item[(A3)] there exist constants $c_A>0$ and $p\geq2$ such that 
  \begin{align*}
   % &\ip{A(\omega,t,\varphi_1)-A(\omega,t,\psi)}{\varphi-\psi}_V\geq0\,,\\
    &\ip{A(\omega,t,\varphi)}{\varphi}_V\geq c_A\norm{\varphi}_V^p\,,
  \end{align*}
  for every $(\omega,t)\in\Omega\times[0,T]$ and $\varphi\in V$;
  \item[(A4)] there exists a constant $C_A>0$
  and a progressively measurable process $h\in L^{1}(\Omega\times(0,T))$ such that,
  setting $q:={p}/{(p-1)}$,
  \[
  \norm{A(\omega,t,\varphi)}^q_{V^*}\leq C_A\norm{\varphi}_V^{p} + h(\omega,t)
  \]
  for every $(\omega,t)\in \Omega\times[0,T]$ and $\varphi\in V$;
  \item[(A5)] for every increasing Lipschitz-continuous function $\sigma\in C^2(\erre)$
  with $\sigma(0)=0$ it holds
  \begin{align*}
  &\sigma(\varphi) \in V \quad\forall\,\varphi\in V\,,\\
  &\sigma_{|V}:V\to V \text{ is locally bounded}\,,\\
  &\ip{A(\omega,t,\varphi) - A(\omega,t,\psi)}{\sigma(\varphi-\psi)}_V\geq 0
  \qquad\forall\,(\omega,t)\in\Omega\times[0,T]\,, \quad\forall\,\varphi,\psi\in V\,.
  \end{align*}
\end{itemize}
 Note that (A1)-(A5) hold, for instance, with the choice $a(\omega,t,\xi) =
\alpha (\omega,t)|\xi|^{p-2}\xi$ with $\alpha$ measurable, bounded, and
uniformly positive.  
 
By choosing  $\sigma(r)=r$  in condition (A5) one in particular has that $A$ is
{\it monotone}. On the other hand, the choice $\sigma(r)=r^+=\max\{r,0\}$
corresponds to the so-called {\it $T$-monotonicity} of $A$, see
\cite{brezis-stampacchia68,Calvert70}. These two functions, together with some
locally regularised version of $r^+$, see \eqref{eq:sigmae}, are
actually the only ones used in the analysis. This would give the
possibility of weakening
assumption (A5), by explicitly referring to these.  

Starting from the 
Carath\'eodory function $b:
\erre \to \erre$, nondecreasing and linearly bounded in the third variable, we define the operator $B :\Omega \times [0,T]\times
H \to H $ as
$$B(\omega,t,u)(x) = b(\omega,t,u(x)) \quad  \text{for a.e.} \
(\omega,t,x)\in \Omega\times[0,T] \times \OO\,.$$
In particular, we require $B$ to fulfill 
\begin{itemize}
\item[(B1)] $B$ is $\mathcal P\otimes\cB(H)$--$\cB(H)$ measurable;
\item[(B2)] $u_1, u_2 \in H, \  u_1 \leq u_2$
  a.e. $\Rightarrow$ $B(\cdot,u_1) \leq B(\cdot,u_2)$ a.e.
\item[(B3)] there exists a constant $C_B>0$ such that
  \[|B(\omega,t,u(x))| \leq C_B\left(1 + | u (x)| \right) \qquad\forall\,
  u \in H, \ \text{for a.e.} \ (\omega,t,x)\in \Omega \times
  (0,T)\times \OO\,.
  \]
\end{itemize} 
 Note that no continuity is required on $b$ nor or $B$.

Again by possibly allowing additional dependencies, we let the
  operator $F:\Omega\times[0,T]\times H \to H$ be defined by 
$$ F(\omega,t,u)(x) = f(\omega,t,u(x)) \quad   \text{for a.e.} \
(\omega,t,x)\in 
\Omega\times[0,T] \times \OO$$
where $f: \Omega \times [0,T]\times \Rz \to \Rz$  is a
Carath\'eodory function, Lipschitz continuous with respect to the last
variable. Specifically, we directly assume on the operator $F$ the following: 
\begin{itemize}
  \item[(F1)] $F$ is $\mathcal P\otimes\cB(H)$--$\cB(H)$ measurable;
  \item[(F2)] there exists a constant $C_F>0$ such that,   
  $$\forall \,u_1,\, u_2 \in H: \quad \|F(\cdot,u_1)-F(\cdot,u_2)\|_H\leq
  C_F \| u_1 - u_2\|_H \quad \text{a.e. in} \ \ \Omega\times (0,T)\,;$$
  \item[(F3)] there exists $\varphi_F\in H$ such that
    $F(\cdot,\cdot,\varphi_F)\in L^2_{\mathcal P}(\Omega; L^2(0,T; H))$.
\end{itemize}

Eventually, the operator $G:\Omega\times[0,T]\times H \to \cL^2(U,H)$  is
required to satisfy 
\begin{itemize}
  \item[(G1)] $G$ is $\mathcal P\otimes\cB(H)$--$\cB(\cL^2(U,H))$ measurable;
  \item[(G2)] there exists a constant $C_G>0$ such that,
  for every measurable subset $\bar\OO\subset\OO$,
  \[
  \sum_{k=0}^\infty\int_{\bar\OO}
  |G(\omega,t,\varphi)e_k-G(\omega,t,\psi)e_k|^2\leq C_G^2\int_{\bar\OO}|\varphi-\psi|^2
  \]
  for every $(\omega,t)\in\Omega\times[0,T]$ and $\varphi,\psi\in H$;
  \item[(G3)] there exists $\varphi_G\in H$ such that $G(\cdot,\cdot,\varphi_G)\in 
  L^2_{\mathcal P}(\Omega; L^2(0,T; \cL^2(U,H)))$.
\end{itemize} 
Assumption (G2) is a generalized Lipschitz-continuity requirement on  $G$.
It is not difficult to check that it is satisfied when  $G$
has the form 
  \[
    G(\omega,t,\varphi)e_k=g_k(\omega,t,\varphi)\,, \quad k\in\enne\,,
  \]
where 
$g_k:\Omega\times[0,T]\times\erre\to\erre$, $k\in\enne$, is Carath\'eodory and 
  \begin{align*} 
  \sum_{k=0}^\infty|g_k(\omega,t,r)-g_k(\omega,t,s)|^2\leq C_G^2|r-s|^2\,,
  \end{align*}
for every $(\omega,t)\in\Omega\times[0,T]$ and $r,s\in\erre$. 

Given the above positions, the variational formulation of (some
extension to additional dependencies of) equation \eqref{eq:i0} along
with variationally defined boundary conditions and an initial condition reads
\beq\label{eq:SPDE!}
  \begin{cases}
  \d u + A(u)\,\d t - B(u)\, \d t= F(u)\,\d t + G(u)\,\d W \quad
  \text{in} \ 
  V^*, \ \text{a.e.~in} \ \Omega \times (0,T) \,,\\
  u(0)=u_0\,,
  \end{cases}
\eeq
As the nonlinear term $-B(u)$ is not monotone and not Lipschitz
continuous, existence for \eqref{eq:SPDE!} does not follow from the
classical theory \cite{KR-spde,LiuRo,Pard}.  
In order to state our existence result, let us first recall a
classical statement on well-posedness in case $B=0$.

\begin{lem}[Case $B=0$]\label{lem:classic}
 Assume  {\em(S1)--(S2)}, {\em(A1)--(A5)}, {\em(F1)--(F3)}, and
 {\em(G1)--(G3)}. For any initial datum $u_0\in L^2(\Omega,\cF_0; H)$ and any $h \in
L^2_{\mathcal P}(\Omega;L^1(0,T;H))$ the Cauchy
problem 
\begin{equation}\label{eq:SPDE}
  \begin{cases}
  \d u + A(u)\,\d t =  h\,\d t + F(u)\,\d t + G(u)\,\d W \quad 
  \text{in} \ 
  V^*, \ \text{a.e. in} \ \Omega \times (0,T) \,,\\
  u(0)=u_0\,
  \end{cases}
\end{equation}
admits a unique solution $u \in L^2(\Omega; C^0([0,T]; H))\cap L^p_{\mathcal P}(\Omega; L^p(0,T; V))$,
in the sense that 
\[
  u(t) + \int_0^t A(s, u(s))\,\d s = u_0 + \int_0^t h(s)\,\d s
  +\int_0^tF(s,u(s))\,\d s + \int_0^tG(s,u(s))\,\d W(s) \quad
  \text{in} \ 
  V^*
\]
for every $t\in[0,T]$, $\P$-almost surely. 
\end{lem}

The crucial tool in our analysis is a comparison principle for
solutions to the Cauchy problem \eqref{eq:SPDE} with respect to the data. We
have the following.

\begin{proposition}[Comparison principle]
  \label{thm:1}
  Assume  {\em(S1)--(S2)}, {\em(A1)--(A5)}, {\em(F1)--(F3)}, and {\em(G1)--(G3)}.
  Let 
  \[
  u_0^1, u_0^2 \in L^2(\Omega, \cF_0; H)\,, \qquad
  h_1,h_2 \in L^2_{\mathcal P}(\Omega; L^1(0,T; H))\,,
  \]
  and let
  \[
  u_1, u_2 \in L^2(\Omega; C^0([0,T]; H))\cap L^p_{\mathcal P}(\Omega; L^p(0,T; V))
  \] 
  be the unique solutions to the Cauchy problem \eqref{eq:SPDE} with respect to 
  data $(u_0^1, h_1)$ and $(u_0^2, h_2)$, respectively. 
  If 
  \begin{align*}
  u_0^1 \leq u_0^2 \quad\text{a.e.~in } \Omega\times\OO\,, \qquad
  h_1\leq h_2 \quad\text{a.e.~in } \Omega\times(0,T)\times\OO\,,
  \end{align*}
  then 
  \[
  u_1(t)\leq u_2(t) \quad\text{a.e.~in } \Omega\times\OO\,, \quad\forall\,t\in[0,T]\,.
  \]
\end{proposition}

The proof of the comparison principle is given in Section
\ref{sec:thm1} and corresponds to an extension of the  
 former analogous result by 
by {\sc Chekroun, Park,\& Temam}
\cite{chek-park-tem} to the case of a nonlinear operator $A$.

As the functions $r \mapsto \pm C_B(1+|r|)$, $r\in \erre$ are
Lipschitz-continuous, owing to Lemma \ref{lem:classic} we can uniquely find
\[
  u_*,u^* \in L^2(\Omega; C^0([0,T]; H))\cap L^p_{\mathcal P}(\Omega; L^p(0,T;V))
\]
solving the Cauchy problems
\begin{align}
  &\begin{cases}
  \d u_* + A(u_*)\,\d t  \\
\hspace{12mm}= - C_B(1+|u_*|)\,\d t +  F(u_*)\,\d t +
  G(u_*)\,\d W\quad \text{in} \ V^*, \ \text{a.e. in} \ \Omega \times (0,T) \,,\\
  u_*(0)=u_0\,,
  \end{cases}\label{eq:ustar}\\[3mm]
  &\begin{cases}
  \d u^* + A(u^*)\,\d t   \\
\hspace{12mm}=  C_B
(1+|u^*|)\,\d t +  F(u^*)\,\d t + G(u^*)\,\d W\quad \text{in} \ V^*, \ \text{a.e. in} \ \Omega \times (0,T) \,,\\
  u_*(0)=u_0\,, 
  \end{cases}\label{eq:usstar}
\end{align}
respectively. Since $-C_B(1+|r|) \leq 0 \leq C_B(1+|r|)$, an
application of Proposition \ref{thm:1} ensures that $u_* \leq u^*$
almost everywhere. 

We can now state our main result on existence of
solutions for the Cauchy problem  \eqref{eq:SPDE!}.

\begin{thm}[Existence]\label{thm:2}
  Assume {\em(S1)--(S2)}, {\em(A1)--(A5)},  {\em (B1)--(B3)}, 
  {\em(F1)--(F3)}, and
  {\em(G1)--(G3)},  
  Then, for any initial datum
  $u_0\in L^2(\Omega,\cF_0; H)$ the Cauchy problem
  \eqref{eq:SPDE!} admits a solution 
  $u \in L^2(\Omega;
  C^0([0,T]; H))\cap L^p_{\mathcal P}(\Omega; L^p(0,T; V))$, in the
  sense that
  \[
  u(t) + \int_0^tA(s, u(s))\,\d s - \int_0^tB(s, u(s))\,\d s = 
  u_0 + \int_0^t F(s,u(s))\,\d s + \int_0^t G(s, u(s))\,\d W(s)
  \quad \text{in} \ V^*
  \]
  for every $t\in[0,T]$, $\P$-almost surely. Moreover, one can uniquely
  find a minimal solution $u_{\rm min}$ and a maximal solution
  $u_{\rm max}$ such that every solution $u$ fulfills $u_* \leq
  u_{\rm min} \leq u \leq u_{\rm max}\leq u^*$ a.e.
\end{thm}

The proof of Theorem \ref{thm:2} is presented in Section
\ref{sec:thm2} and relies on a fixed-point procedure for nondecreasing
mappings. Note that no uniqueness for the Cauchy problem
\eqref{eq:SPDE!}  is to be expected. Indeed, the classical
counterexample to uniqueness in $\erre$ given by the deterministic ODE problem
$$u'=(\max\{u,0\})^{1/2}, \quad u_0=0,$$ is included in the setting of
Theorem \ref{thm:2}. In this case, $u_{\rm min}(t)=0$ and $u_{\rm max}(t)
= t^2/4$ for $t \geq 0$.

%%%%%%%%%%%%%%%%%%%%%%%%%%%%%%%%%%%%%%%%%%%%%%%%%%

\section{Comparison principle: proof of Proposition \ref{thm:1}}

\setcounter{equation}{0}
\label{sec:thm1}

We closely follow here the argument from \cite{chek-park-tem}, by
adapting it to our nonlinear setting. Under the notation of
Proposition~\ref{thm:1}, we introduce the new variable $u:=u_1-u_2$
and define
$h:=h_1-h_2$ and $u_0^:=u_0^1-u_0^2$. Then, $u$ satisfies the Cauchy problem
\[
  \begin{cases}
  \d u + (A(u_1) - A(u_2))\,\d t \\
\qquad = h\,\d t + (F(u_1)-F(u_2))\,\d t
  +(G(u_1)-G(u_2))\,\d W \quad \text{in $V^*$, a.e. in
    $\Omega\times (0,T)$}\,,\\
  u(0) = u_0\,.
  \end{cases}
\]

Introduce now the operators 
\begin{align*}
  &\tilde F:\Omega\times[0,T]\times H \to H\,, \\
  &\tilde F(\omega,t,\varphi):=F(\omega,t,\varphi + u_2(\omega,t)) - F(\omega,t, u_2(\omega,t))\,,
  \quad(\omega,t,\varphi) \in \Omega\times[0,T]\times H\,,\\[3mm]
  &\tilde G:\Omega\times[0,T]\times H \to \cL^2(U,H)\,, \\
  &\tilde G(\omega,t,\varphi):=G(\omega,t,\varphi + u_2(\omega,t)) - G(\omega,t, u_2(\omega,t))\,,
  \quad(\omega,t,\varphi) \in \Omega\times[0,T]\times H\,.
\end{align*}
Note that 
$\tilde F$ and $\tilde G$ still satisfy
assumptions (F1)--(F3) and (G1)--(G3), respectively.
Additionally, by definition we have
\beq
  \label{FG_0}
  \tilde F(\cdot,\cdot,0)=0\,, \qquad
  \tilde G(\cdot, \cdot, 0)=0\,.
\eeq
With this notation, the Cauchy problem for $u$ can be
equivalently rewritten
as
\beq\label{cau_diff}
  \begin{cases}
  \d u + (A(u_1) - A(u_2))\,\d t = h\,\d t + \tilde F(u)\,\d t
  +\tilde G(u)\,\d W  \quad \text{in $V^*$, a.e. in $(0,T)$, \
  $\P$-a.s.}\\
  u(0) = u_0\,.
  \end{cases}
\eeq

Recall that we have $u_0\leq 0 $ a.e.~in $ \Omega\times\OO$ and 
$h\leq0 $ a.e.~in $\Omega\times(0,T)\times\OO$.
Along with this notation, the assertion follows by proving that
$u(t)\leq 0$ a.e.~in $\Omega\times\OO$ for all $t\in[0,T]$.
We check this by showing that
\beq
  \label{sigma_0}
  (u(t))^+=0
  \qquad\text{a.e.~in } \Omega\times\OO\,, \quad\forall\,t\in[0,T]\,.
\eeq

In order to prove \eqref{sigma_0}, we resort in an approximation
of the positive part by means of the sequence
$(\sigma_{\eps})_{\eps>0}$,
defined in \cite[\S~2.4]{chek-park-tem} as  
\begin{equation}
  \sigma_\eps(r):=\begin{cases}
  r \quad&\text{if } u>\eps\,,\\[1mm]
  \displaystyle \frac3{\eps^4} r^5 - \frac8{\eps^3}r^4 + \frac6{\eps^2}r^3 \quad&\text{if } 0<r\leq\eps\,,\\[2mm]
  0 \quad&\text{if } r<0\,.
  \end{cases}
\label{eq:sigmae}
\end{equation}
It is not difficult to check that $\sigma_\eps\in C^2(\erre)$ for every $\eps>0$,
and that there exists a constant $M>0$, independent of $\eps$, such that 
\beq
  \label{bound_sig}
  |\sigma_\eps'(r)| + |\sigma_\eps''(r)| + |\sigma_\eps(r)\sigma_\eps''(r)| \leq M 
  \qquad\forall\,r\in\erre\,, \quad\forall\,\eps>0\,.
\eeq
Moreover, $\sigma_\eps\geq 0$ for every $\eps>0$ and $\sigma_\eps(r)\nearrow r^+$
for all $r\in \erre$ as $\eps\searrow0$.
Defining now the primitive functions
\[
  \hat\sigma_\eps:\erre\to\mathopen[0,\infty)\,, \qquad
  \hat\sigma_\eps(r):=\int_0^{r}\sigma_\eps(s)\,\d s\,, \quad r\in\erre\,,
\]
we introduce the functional $\Sigma_\eps: H \to
\mathopen[0,\infty\mathclose)$ as 
\[ 
  \Sigma_\eps(\varphi):=\int_\OO \hat\sigma_\eps(\varphi)\,, \quad \varphi\in H\,.
\]
We aim now at applying It\^o's formula to $\Sigma_\eps(u)$.
This is indeed possible since $\Sigma_\eps$ is Fr\'echet
differentiable in $H$, with 
derivative given by 
\[
  D\Sigma_\eps: H\to H\,, \qquad D\Sigma_\eps(\varphi)=\sigma_\eps(\varphi)\,,
  \quad\varphi\in H\,.
\]
Moreover, since $V\embed L^4(\OO)$, it follows that the restriction of 
$D\Sigma_\eps$ to $V$ is Fr\'echet differentiable in $V$ and its derivative is given by
\[
  D^2\Sigma_\eps:V\to \cL(V,H)\,, \qquad
  D^2\Sigma_\eps(\varphi)w=\sigma_\eps'(\varphi)w\,,
  \quad\varphi,w\in V\,.
\]
From (A5) we have that the 
restriction of $D\Sigma_\eps$ to $V$ takes values in $V$, and that 
$D{\Sigma_\eps}_{|V}:V\to V$ is strongly-weakly continuous.
We can hence apply It\^o's formula in the variational
setting of \cite[Thm.~4.2]{Pard} and obtain
\begin{align*}
  &\Sigma_\eps(u(t)) + \int_0^t\ip{A(s,u_1(s))-A(s,u_2(s))}{\sigma_\eps(u(s))}_V\,\d s \\
  &\qquad=\Sigma_\eps(u_0) + \int_0^t\left(h(s), \sigma_\eps(u(s))\right)_H\,\d s
  +\int_0^t\left(\tilde F(s,u(s)), \sigma_\eps(u(s))\right)_H\,\d s\\
  &\qquad+\int_0^t\left(\sigma_\eps(u(s)), \tilde
    G(s,u(s))\,\d W(s)\right)_H
  +\frac12\int_0^t\sum_{k=0}^\infty
  \left(\sigma_\eps'(u(s))\tilde G(s,u(s))e_k, \tilde G(s,u(s))e_k\right)_H\,\d s
\end{align*}
for every $t\in[0,T]$, $\P$-almost surely.
Since $u_0\leq 0$ almost everywhere we have 
$\Sigma_\eps(u_0)=0$. Moreover, since 
$h\leq0$ and $\sigma_\eps(u)\geq0$ almost everywhere, 
the second term on the right-hand side is nonpositive.
Noting also that the second term on the left-hand side is
nonnegative by (A5), by
taking expectations we infer that 
\begin{align*}
  \E\Sigma_\eps(u(t))
  &\leq
  \E\int_0^t\left(\tilde F(s,u(s)), \sigma_\eps(u(s))\right)_H\,\d s\\
  &+\frac12\E\int_0^t\sum_{k=0}^\infty
  \left(\sigma_\eps'(u(s))\tilde G(s,u(s))e_k, \tilde G(s,u(s))e_k\right)_H\,\d s\,.
\end{align*}
Now, by definition of $\sigma_\eps$,
the uniform estimates \eqref{bound_sig}, 
assumptions (F3) and (G3), and the Dominated Convergence Theorem,
letting $\eps\searrow0$ we infer that 
\begin{align*}
 \frac12\E\|u^+(t)\|_H^2 &\leq 
 \E\int_0^t\left(\tilde F(s,u(s)), u^+(s) \right)_H\,\d s
  +\frac12\E\int_0^t\sum_{k=0}^\infty\int_{\{u(s)\geq 0\}}|\tilde G(s,u(s))e_k|^2\,\d s\\
  &=\E\int_0^t\!\!\int_{\{u(s)\geq0\}}\left(\tilde
    F(s,u(s)),u^+(s)\right)_H \,\d s
  +\frac12\E\int_0^t\sum_{k=0}^\infty\int_{\{u(s)\geq 0\}}|\tilde G(s,u(s))e_k|^2\,\d s\,
\end{align*}
for all $t\in [0,T]$.
By using the H\"older inequality, 
the Lipschitz-continuity assumptions (F2) and (G2) on $\tilde F$ and $\tilde G$, 
together with the fact that 
$\tilde F(\cdot,0)=\tilde G(\cdot, 0)=0$ from \eqref{FG_0}, we deduce that 
\begin{align*}
\frac12\E\|u^+(t)\|_H^2 &\leq 
  % C_F\E\int_0^t\left(\int_{\{u(s)\geq0\}}|u(s)|^2\right)^{1/2}\norm{u^+(s)}_H\,\d s
  % +C_G^2\E\int_0^t\!\!\int_{\{u(s)\geq 0\}}|u(s)|^2\,\d s\\
  (C_F + C_G^2)\int_0^t \E\|u^+(s)\|_H^2\,\d s\,.
\end{align*}
Hence, \eqref{sigma_0} follows from the Gronwall lemma, and Theorem~\ref{thm:1} is proved.

%%%%%%%%%%%%%%%%%%%%%%%%%%%%%%%%%%%%%%%%%%%%%%%%%%

\section{Existence of solutions: proof of Theorem \ref{thm:2}}

\setcounter{equation}{0}
\label{sec:thm2} 

As anticipated, the proof of Theorem \eqref{thm:2} relies on 
a fixed-point tool for nondecreasing mappings in ordered sets. Let us
start by recalling some basic notion.

Let $(E, \preceq) $ denote a nonempty ordered set and $ F\subset E$. We
recall that $ f \in F $ is a {\it maximal (minimal) element} of $ F$
iff, for all $f'\in F$, $ f\preceq f'$ ($ f'\preceq f$, respectively)
implies $f=f'$ and that $ f $ is the {\it maximum (minimum)} of $ F $
iff $ f'\preceq f$ ($ f\preceq f'$, respectively) for all $ f' \in
F$. Moreover, $ e \in E$ is an {\it upper bound (lower bound)} of $ F
$ iff $ f\preceq e$ ($ e\preceq f$, respectively) for all $ f \in F$ and $
e\in E $ is the {\it supremum} or {\it least upper bound} ({\it
  infimum} or {\it greatest lower bound}) iff $ e $ is the minimum
(maximum) of the set of upper bounds (lower bounds, respectively) of
$F$. Eventually, we say that $F $ is a {\it chain} if it is totally
ordered and that $ F $ is an {\it interval} iff there exist $ e_*,e^*
\in E $ such that $ F\equiv \{ e \in E : e_* \preceq e \preceq e^*\}$. In
the latter case we use the notation $ F=[e_*,e^*]$. The set $(E, \preceq)
$ is said to be {\it s-inductive (i-inductive)} iff every chain of $ E
$ is bounded above (below, respectively) and $ (E, \preceq) $ is said to
be {\it completely s-inductive (completely i-inductive)} iff every
chain of $ E $ has a supremum (infimum, respectively). Finally $ (E,
\preceq) $ is said to be {\it inductive (completely inductive)} iff it is
both {s-inductive} and {i-inductive} ({completely s-inductive} and 
{completely i-inductive}, respectively). 

Let us choose $E:=L^2_{\mathcal P}(\Omega; L^2(0,T;
H))$ and specify 
\[
  v_1\preceq v_2 \quad\text{iff}\quad v_1\leq v_2 \quad\text{a.e.~in } \Omega\times(0,T)\times\OO\,,
  \qquad v_1,v_2\in E\,.
\]
By fixing a
tentative $\tilde u \in E$ in the nonlinearity $-B(\tilde u)$, one
recalls assumptions (B1) and (B3) 
giving $B(\tilde u)\in L^2_{\mathcal P}(\Omega;L^2(0,T;H))$. 
By using Lemma \ref{lem:classic}, one uniquely finds 
$$u \in L^2(\Omega; C^0([0,T]; H))\cap L^p_{\mathcal P}(\Omega; L^p(0,T; V))\subset E$$ 
solving the Cauchy problem  
$$\begin{cases}
  \d u + A(u)\,\d t =  B(\tilde u)\,\d t + F(u)\,\d t + G(u)\,\d W \quad
  \text{in} \ 
  V^*, \ \text{a.e. in} \ \Omega \times (0,T) \,,\\
  u(0)=u_0\,.
  \end{cases}
$$
This defines a mapping $S:E \to E$ as 
$$S(\tilde u):=u\,.$$
 The function $u \in E$ is hence a solution of the Cauchy problem 
\eqref{eq:SPDE!} if and only if it is a fixed point of $S$. We will use the following fixed-point lemma.

\begin{lem}[Fixed point]\label{kolo} Let $ (E, \preceq) $ be an ordered set and 
  $I:= [e_*,e^*] \subset E $ be completely inductive. Suppose that
  $S: (I,\preceq) \rightarrow (I,\preceq)\,$ is nondecreasing. Then, the
  set of fixed points $ \{ u \in I \ : \ u = S(u)\} $ is nonempty and has a
  minimum and a maximum.
\end{lem}

This fixed-point result was announced by {\sc Kolodner} \cite{Kolodner68} and turns out to be the main tool in the analysis of \cite{Mignot-Puel77,Tartar74}. 
Its proof is to be found, for instance, in \cite[Thm. 9.26, p. 223]{Baiocchi-Capelo84}.
This fixed-point lemma corresponds indeed to an
abstract version of the classical Perron's method. In particular, in
order to identify the unique minimal fixed point  of $S$ one subsequently
proves that  the set of {\it subsolutions} $ A:=\{v \in I \ : \ v
\preceq S(v) \}$ is non-empty,  $A$ with the induced order is
completely s-inductive,   $A$ has a maximal element $ u$, and $u$ is a fixed point for $ S$.

In order to apply the fixed-point Lemma \ref{kolo} we define $e_*=u_*$
and $e^*=u^*$, where $u_*$ and $u^*$ are the unique solutions to
\eqref{eq:ustar} and \eqref{eq:usstar}, respectively,  and check
that (1) $I$ is completely inductive, (2) $S$ is nondecreasing, and (3)
$S(I)\subset I$.  

Ad (1): Let
$\emptyset \not = F
\subset I$ be a chain. For almost all $(\omega,t
,x)\in \Omega \times (0,T)\times \OO$ we have $(\sup F)(\omega,t ,x) =\sup\{u (\omega,t
,x) \ | \ u \in F\} $ and $(\inf F)(\omega,t ,x) =\inf\{u (\omega,t
,x) \ | \ u \in F \}$, so that $\sup F,\, \inf F \in I$. Hence, $I$ is
completely inductive. 

Ad (2): Take $\tilde
u_1 \preceq \tilde u_2$ and recall that $u_1 = S(\tilde u_1)$ and $u_2
= S(\tilde u_2)$ are the unique solutions to the Cauchy problem
\eqref{eq:SPDE} with $h$ replaced by $ h_1=B(\tilde u_1)$ and $h_2=B(\tilde
u_2)$, respectively. As $B$ is nondecreasing, we have that $h_1
\preceq h_2$. By applying Proposition \ref{thm:1} we then find  $u_1
\preceq u_2$. This proves that   $S(\tilde u_1)\preceq S(\tilde u_2)$,
namely $S$ is nondecreasing.

Ad (3): Let $\tilde u\in I$
and set $u = S(\tilde u)$.  As  $u_*\preceq \tilde u$ and $B$ is
nondecreasing, we have that $B(u_*)\preceq B(\tilde u)$. Assumption (B2) ensures that  we have that 
 \beq\label{alpha_bound}
   -C_B(1+|v|) \leq |B(\cdot,v)|\leq C_B(1+|v|)
   \quad\forall\,v \in H, \ \text{a.e. in} \ \Omega \times (0,T)\times \OO\,.
 \eeq
Consequently,  we deduce that 
\[
  -C_{B}(1+|u_*|)\leq B(\cdot, u_*)\leq B(\cdot,
  \tilde u) \quad\text{a.e.~in } \Omega\times(0,T)\times\OO\,.
\]
Noting also that $u_*(0)= \tilde u (0)=u_0$, we can apply
Proposition~\ref{thm:1} with the choices
$u_0^1=u_0^2=u_0$, $h_1=-C_{B}(1+|u_*|)$, and $h_2=
B(\cdot,\tilde u)$ and deduce that  $u_*\preceq S(\tilde u)$. An analogous argument
entails the upper bound $S(\tilde u) \preceq u^*$, so that $u_*
\preceq S(\tilde u) \preceq u^*$ or, equivalently, $S(\tilde u)\in I$.

We are hence in the position of applying Lemma~\ref{kolo} and
find that fixed points of $S$ in $I$ exists and admit (unique) maximum
and minimum. The proof
of Theorem~\ref{thm:2} follows then by checking that {\it all}
solutions $u$ to the Cauchy problem \eqref{eq:SPDE!} necessarily
belong to $I$. This follows by applying once again Proposition
\ref{thm:1} and using relations \eqref{alpha_bound}.  
 
%%%%%%%%%%%%%%%%%%%%%%%%%%%%%%%%%%%%%%%%%%%%%%%%%%%%%

\section*{Acknowledgement}
This work has been funded by the Vienna Science and Technology Fund (WWTF)
through Project MA14-009.
The support by the Austrian Science Fund (FWF) projects 
F\,65 and I\,2375 is also gratefully acknowledged. 

%%%%%%%%%%%%%%%%%%%%%%%%%%%%%%%%%%%%%%%%%%%%%%%%%%%%%

\bibliographystyle{abbrv}

\def\cprime{$'$}

\end{document}